\documentclass[10pt]{article}

\usepackage{amsmath,amssymb,mathtools,amsthm}
\usepackage{graphicx}
\usepackage{xcolor}

\newcommand {\eq} [1] {\begin{equation}\label{#1}}
\newcommand {\en} {\end{equation}}
\newcommand {\cB}       {{\cal B}}

\newcommand {\cE}       {{\cal E}}

\newcommand {\cR}       {{\cal R}}
\newcommand {\cS}       {{\cal S}}

\newcommand {\cX}       {{\cal X}}

\newcommand {\eproof}
      {\space
        {\ \vbox{\hrule\hbox{\vrule height1.3ex\hskip0.8ex\vrule}\hrule}}
        \par}
\newcommand {\C}        {{\mathbb C}}
\newcommand {\R}        {{\mathbb R}}

\newcommand {\Rnn}      {\R^{n \times n}}
\newcommand {\Rnm}      {\R^{n \times m}}
\newcommand {\Rmm}      {\R^{m \times m}}
\newcommand {\Rmn}      {\R^{m \times n}}

\newcommand {\mat}      [1] {\left[\begin{array}{#1}}
\newcommand {\rix}          {\end{array}\right]}
\newtheorem{theorem}           {Theorem}
\newtheorem{lemma}    [theorem]{Lemma}

\newcommand {\rank}     {\mathop{\rm rank}\nolimits}

 \catcode`@=11
 \font\tenex=cmex10 
 \newdimen\p@renwd
 \setbox0=\hbox{\tenex B} \p@renwd=\wd0 
 \def\bmat#1{\begingroup \m@th
   \setbox\z@\vbox{\def\cr{\crcr\noalign{\kern2\p@\global\let\cr\endline}}%
     \ialign{$##$\hfil\kern2\p@\kern\p@renwd&\thinspace\hfil$##$\hfil
       &&\quad\hfil$##$\hfil\crcr
       \omit\strut\hfil\crcr\noalign{\kern-\baselineskip}%
       #1\crcr\omit\strut\cr}}%
   \setbox\tw@\vbox{\unvcopy\z@\global\setbox\@ne\lastbox}%
   \setbox\tw@\hbox{\unhbox\@ne\unskip\global\setbox\@ne\lastbox}%
   \setbox\tw@\hbox{$\kern\wd\@ne\kern-\p@renwd\left[\kern-\wd\@ne
     \global\setbox\@ne\vbox{\box\@ne\kern2\p@}%
     \vcenter{\kern-\ht\@ne\unvbox\z@\kern-\baselineskip}\,\right]$}%
   \null\;\vbox{\kern\ht\@ne\box\tw@}\endgroup}
 \catcode`@=12
\begin{document}

\title{Asymptotic Stability and Strict Passivity of Port-Hamiltonian Descriptor Systems \\ via State Feedback}
\author{Delin Chu\footnotemark[1]
\and Volker Mehrmann\footnotemark[2]}

\maketitle

\renewcommand{\thefootnote}{\fnsymbol{footnote}}
\footnotetext[1]{Department of Mathematics, National University of Singapore, Singapore 119076. Email: {\tt
 matchudl@nus.edu.sg}.}
\footnotetext[2]{Institute of Mathematics, TU Berlin, Stra{\ss}e des 17.~Juni 136, 10623 Berlin, Germany.  Email: {\tt mehrmann@math.tu-berlin.de}.}

\begin{abstract} While port-Hamiltonian descriptor systems are known to be stable and passive, they may not be asymptotically stable or strictly passive. Necessary and sufficient conditions  are presented when these properties  as well as the regularity and the index one property can be achieved via state feedback while preserving the port-Hamiltonian structure.
\end{abstract}

{\bf Keywords:}   Port-Hamiltonian descriptor system, proportional state feedback, regularity, index reduction, asymptotic stability, strict passivity.

 {\bf AMS subject classification:} 93B05, 93B40, 93B52, 65F35

\section{Introduction}\label{intro}

In this paper we consider linear time-invariant port-Hamiltonian (pH) descriptor systems of the form
\begin{eqnarray}
E {\dot x} &=& (J-R)x+(G-P) u,  \label{genphdae}  \\
         y &=&  (G+P)^Tx+(S+N)u,    \nonumber
\end{eqnarray}
where $ E, J, R\in \Rnn$, $G, P\in \Rnm$,  $D=S+N\in \Rmm$ ($S=S^T=\frac 12(D+D^T)$, $N=-N^T=\frac 12(D-D^T)$), and the functions  $x,y,u$ are the state, output,  and input  of the system. The coefficient matrices satisfy
\[
 E = E^T\geq 0, \ \mat{cc} J & G \\ -G^T & N \rix=-\mat{cc} J & G \\ -G^T & N \rix^T,   \]
\begin{equation} W= W^T=\mat{cc} R & P \\ P^T & S \rix\geq 0, \label{lmi-1}
\end{equation}
where for a matrix $M$, $M \geq 0$ $(M>0)$ denotes that $M$ is positive semi-definite (definite) and $\Rnm$ denotes the real $n\times m$ matrices.

{Note that there are different representations of pH descriptor systems, some include an extra term $Q$ in front of the state $x$ in \eqref{genphdae}. It has been discussed in
several publications, see \cite{MehMW18,MehU23,MehS23} how the system can be transformed to remove this term, so in the following we use the form \eqref{genphdae}.
}

The class of pH descriptor  systems represents an ideal  modeling framework for control systems in almost all real world physical systems. { While ordinary pH systems are well studied see e.~g.~\cite{JacZ12,RasCSS20,SchJ14}, the study of pH descriptor systems has only recently received a lot of attention, see~\cite{BeaMXZ18,MehM19,MehS23,MehU23,Mor24,Sch13,SchM18,SchM20,SchM23} for a detailed discussion of the many important properties and a large number of applications.}
{ These properties include the stability and passivity of the system. The extension of classical analytical definition of stability for ordinary differential equations, see e.~g. \cite{HinP05} to descriptor systems is well-known, see e.~g. \cite{KunM24,MehU23}. The difficulty in this extension is that either the set of perturbations in the initial conditions has to be restricted or the class of solution functions has to be enriched to deal with inconsistent initial conditions, see \cite{KunM24} for a detailed analysis. For the case of linear time-invariant descriptor systems, however, it is a well established fact, see e.~g. \cite{DuLM13}, that the stability analysis can be carried via the eigenstructure of the associated matrix pencil which in our case is the \emph{dissipative Hamiltonian matrix pencil} $sE-(J-R)$.}



A characteristic property of systems of the form \eqref{genphdae} is the availability of a \emph{storage function (Hamiltonian)} which is given by $\mathcal H(x)=\frac 12 x^TEx\geq 0$ and that the following power balance equation holds, as well as the resulting dissipation inequality,
see e.~g. \cite{MehM19}.
\begin{theorem}\label{thm:pbe}
Consider a pH descriptor system of the form \eqref{genphdae}. Then for any input function $u$, the \emph{power balance equation}
  \begin{equation}\label{eq:powerBalanceEq}
 {   \frac{d}{dt}\mathcal H(x)     }
            = - \begin{bmatrix}
        x\\ u
    \end{bmatrix}^T\mat{cc} R & P \\ P^T & S \rix\begin{bmatrix}
        x\\ u
    \end{bmatrix}+ y^Tu
  \end{equation}
  holds along any solution $x$.
  In particular, the \emph{dissipation inequality}
  \begin{equation}\label{eq:dissIneq}
    \mathcal H(x(t_2)) - \mathcal H(x(t_1)) \leq \int_{t_1}^{t_2} (y(\tau)^Tu(\tau))\ d\tau
  \end{equation}
  holds.
\end{theorem}
%
%
%
The dissipation inequality \eqref{eq:dissIneq} immediately implies that the system is \emph{passive}, see \cite{CheGH23,Wil71}, but since the inequality in \eqref{eq:dissIneq} may not be strict, the system may not be \emph{strictly passive}.

Similarly, it is well known, see e.~g. \cite{MehMW18}, that  pH descriptor systems of the form \eqref{genphdae} are \emph{stable}, in the sense that all finite eigenvalues of $sE-(J-R)$ are in the closed left half plane and those on the imaginary axis are semi-simple.  But since purely imaginary eigenvalues may occur, pH descriptor systems are not necessarily \emph{asymptotically stable}.

Furthermore, it is shown in  \cite{MehMW18} that the \emph{index} of a pH descriptor system, i.~e. the size of the largest Jordan block associated with the eigenvalue $\infty$ of $sE_(J-R)$, can be at most $2$, and it is shown in \cite{MehMW21} that the \emph{singularity} of the pencil $sE-(J-R)$, i.~e. that
$\det (sE-(J-R))$ is identically $0$, is characterized by a common nullspace  of the coefficients $E,J,R$.

It is well known \cite{Kai80} that neither (asymptotic) stability nor
(strict) passivity and also neither the regularity and the index of a descriptor system are invariant under feedback. On the contrary, typically feedback is used to make an unstable or non-passive system asymptotically stable, respectively strictly passive and to make a descriptor system regular and of index at most one, \cite{ChuMN99}.

{ In \cite{ChuM24a} the difficult question of obtaining a regular, index at most one, and asymptotically stable closed pH system via output feedback has been discussed and a complete characterization has been given. Typically, however, one can achieve a lot better results when one uses state feedback, because for the changes in the system pencil one has a larger class available. On the other hand, state feedback will usually destroy the pH structure.}

The goal of this paper, is therefore, to study the question  whether a pH descriptor system can be made regular, of index at most one, and furthermore strictly passive, respectively asymptotically stable, by applying a state feedback $u=Fx +v$ with $F\in  \Rmn$ to (\ref{genphdae}) while preserving the pH structure.

Applying such a state feedback, the closed-loop system has the form
\begin{eqnarray}
E \dot x &=& [J-R+(G-P)F] x+(G-P)v=(\tilde J-\tilde R) x+(\tilde G-\tilde P) v, \label{clo-sys1}\\
          y &=& (G^T+P^T+(S+N) F)x+(S+N)u=(\tilde G+\tilde P)^T x+(S+N)v,  \nonumber
\end{eqnarray}
where
\begin{eqnarray*}
\tilde J&=&J+{1\over 2}[ (G-P)F-F^T(G-P)^T], \\
   \tilde R&=&R-{1\over 2}[(G-P)F+F^T(G-P)^T], \\
  \tilde G&=&G+{1\over 2} F^T(S+N)^T, \ \tilde P=P+{1\over 2} F^T(S+N)^T.
\end{eqnarray*}
We present a complete characterization for the following two questions.
%
 \begin{itemize}
\item []\label{P1} \textbf{Question~1}:
When does there exist a state feedback matrix $F$ such that  the closed-loop system (\ref{clo-sys1}) associated with \eqref{genphdae} is port-Hamiltonian, regular, of index at most one, and asymptotically  stable.

\item []\label{P2} \textbf{Question~2}: When does there exists a state feedback matrix $F$ such that  the closed-loop system (\ref{clo-sys1}) associated with \eqref{genphdae} is port-Hamiltonian, regular, of index at most one, and strictly passive, i.~e., the dissipation matrix of the closed loop system satisfies
\begin{equation}\label{pa1} \tilde W=\tilde W^T=\mat{cc}  \tilde R & \tilde P \\ \tilde P^T & S \rix>0.
\end{equation}
\end{itemize}

The paper is organized as follows. In Section~\ref{sec:prelim} we present some preliminary results. In Section~\ref{sec:main} we present the answers to the two questions
and we finish with some conclusions and open problems in Section~\ref{sec:conclusions}.

\section{Preliminaries}\label{sec:prelim}
In  this section we introduce some notation and present some preliminary results.

We denote the open left (right) half complex plane by $\C_- (\C_+)$ and the imaginary axis by $\C_0$.
A column orthogonal matrix with its columns spanning the right nullspace of a matrix $M$ is denoted by $\cS_\infty(M)$ and with its columns spanning the range space of $M$ by $\cR(M)$. The Moore-Penrose inverse of a matrix $M$ is denoted by $M^{(+)}$ and the unique positive definite \emph{square root} of a matrix $M=M^T>0$ by $M^{\frac 12}$.

We have the following elementary lemmas.
\begin{lemma}\label{lemma-1} \cite{ChuM24a} Consider $E, J, R\in \R^{n,n}$ with $E=E^T\geq 0$, $J=-J ^H$, $R=R^T\geq 0$, partitioned in blocks as
\[
E=\bmat{ & n_1 & n_2 \cr
             n_1 & E_{11} & E_{12} \cr
             n_2 & E_{12} ^T & E_{22} \cr}, \quad
   R=\bmat{ & n_1 & n_2 \cr
             n_1& R_{11} & 0 \cr
             n_2 & 0         & 0   \cr}, \quad
   J=\bmat{ & n_1 & n_2 \cr
             n_1 & J_{11} & J_{12} \cr
             n_2 & -J_{12} ^T & J_{22} \cr}, \]
where $R_{11}>0$. Then  the following statements hold:

(i) The matrix pencil    $sE-(J-R)$ has all its finite eigenvalues in the open left half complex plane if and only if
\begin{equation}\label{stable-1}
\rank \mat{cc} s E_{12}^T +J_{12}^T &  s E_{22}-  J_{22} \rix=n_2\  \mbox{\textrm for all }\ s\in \C_0.
\end{equation}

(ii) The matrix $J-R$ is nonsingular if and only if
\[
\rank \mat{cc} -J_{12}^T & J_{22} \rix=n_2.
\]
\end{lemma}
A simple { construction} such that a pH descriptor system can be made asymptotically stable via state feedback is given by the following lemma.
\begin{lemma}\label{lemma-2}  Consider $E, J, R\in \Rnn$ and $B\in \Rnm$, with $E=E^T\geq 0$, $R=R^T\geq 0$ and $J=-J^T$, and suppose that
\begin{equation}\label{stable} \rank \mat{cc} sE-(J-R) & B \rix=n\  \mbox{\textrm for all }\ s\in \C_0.
\end{equation}
Let $F\in \Rnm$ be such that $\tilde R:=R-{1\over 2}(BF+F^TB^T)\geq 0$ and
\begin{equation}\label{F1}
\rank \tilde R=\rank \mat{cc} R & B \rix.
\end{equation}
Then the pencil $sE-(J-R+BF)$ has all its eigenvalues in $\C_-$.
\end{lemma}
\proof
{
Note that since $R=R^T$, by first performing an echelon form of $B$ and  then compressing $R$ in the part associated with the left nullspace of $B$, there exist nonsingular matrices $X_1$ and $Y$   such that
\[
X_1RX_1^T=\bmat{ & n_1 & n_2 & n_3 \cr
            n_1 & \cR_{11} & R_{12} & R_{13} \cr
            n_2 & R_{12}^T          & R_{22} &  0 \cr
            n_3 & R_{13}^T & 0 & 0 \cr}, \quad
    X_1BY^T=\bmat{ & n_1 & m-n_1 \cr
              n_1& I     & 0 \cr
              n_2 & 0 & 0 \cr
              n_3 & 0 & 0 \cr},
\]
where $R_{22}$ is nonsingular. Since $R=R^T\geq 0$, it holds that $R_{13}=0$ and $R_{22}>0$.
with
\[ X=\mat{ccc} I & -R_{12}R_{22}^{-1} & 0 \\ 0 & R_{22}^{-1/2} & 0 \\ 0 & 0 & I \rix X_1, \quad
    R_{11}=\cR_{11}-R_{12}R_{22}^{-1}R_{12}^T, \]
we obtain
\[
XRX^T=\bmat{ & n_1 & n_2 & n_3 \cr
            n_1 & R_{11} & 0 & 0 \cr
            n_2 & 0         & I &  0 \cr
            n_3 & 0 & 0 & 0 \cr}, \quad
    XBY^T=\bmat{ & n_1 & m-n_1 \cr
              n_1& I     & 0 \cr
              n_2 & 0 & 0 \cr
              n_3 & 0 & 0 \cr},
\]
and $R_{11}\geq 0$.
}
Transforming and partitioning the other matrices accordingly as
\[
XEX^T=\bmat{ & n_1 & n_2 & n_3 \cr
            n_1 & E_{11} & E_{12} & E_{13} \cr
            n_2 & E_{12}^T & E_{22} & E_{23} \cr
            n_3 & E_{13}^T & E_{23}^T & E_{33} \cr}, \quad
     XJX^T=\bmat{ & n_1 & n_2 & n_3 \cr
             n_1 & J_{11} & J_{12} & J_{13} \cr
             n_2 & -J_{12}^T & J_{22} & J_{23} \cr
             n_3 & -J_{13}^T & -J_{23}^T & J_{33} \cr},
\]
and
\[
Y^{-T}FX^T=\bmat{ & n_1 & n_2 & n_3 \cr
           n_1  & F_{11} & F_{12} & F_{13} \cr
        m-n_1 & F_{21} & F_{22} & F_{23} \cr},
\]
it follows from { the condition $\tilde R:=R-{1\over 2}(BF+F^TB^T)\geq 0$ }
that $F_{13}=0$ and
\[\mat{cc} R_{11}-{1\over 2}(F_{11}+F_{11}^T) & -{1\over 2}F_{12} \\ -{1\over 2}F_{12}^T & I \rix >0.
\]
Then
\[
X \tilde J X^T:=X(J+{1\over 2}(BF-F^TB^T))X^{T}=\mat{ccc} J_{11}+{1\over 2}(F_{11}-F_{11}^T) & J_{12}+{1\over 2}F_{12} & J_{13} \\                 -J_{12}^T-{1\over 2}F_{12}^T & J_{22} & J_{23} \\   -J_{13}^T & -J_{23}^T & J_{33} \rix,
\]
and condition (\ref{stable}) implies that
\[
\rank
\mat{ccc} sE_{13}^T+J_{13}^T & sE_{23}^T+J_{23}^T & sE_{33} -J_{33}\rix=n_3\  \mbox{\textrm for all }\ s\in \C_0.
\]
Therefore, using Lemma \ref{lemma-1} (i) we have that   the pencil $sXEX^T- X(J-R+BF)X^T$ and therefore also the equivalent pencil $sE -(J-R+BF)$ has all its finite eigenvalues in $\C_-$.
\eproof
%

%

{ The following lemma  is important for analyzing when  a pH descriptor system can be made regular and index at most one.
}
 \begin{lemma}\label{lemma-3}  Given $E, J, R\in \Rnn$, and $B\in \Rnm$ with $E=E^T\geq 0$, $R=R^T\geq 0$ and $J=-J^T$. Suppose that
\begin{equation}\label{index-con} \rank \mat{ccc} E & (J-R)\cS_\infty(E) & B \rix=n,
\end{equation}
and let $F\in \Rmn$ be such that $\tilde R:=R-{1\over 2}(BF+F^TB^T)\geq 0$ and
\begin{equation}\label{F2}
    \rank(\tilde R)=\rank \mat{cc} R & B \rix.
\end{equation}
Then the pencil $sE-(J-R+BF)$ is regular and of index at most one.
\end{lemma}
\proof   {
 By first compressing the matrix $B$, next compressing  the matrix $R$ in the left nullspace of $B$ and then compressing
the matrix $E$ in the left  nullspace of $\mat{cc} R & B \rix$, using the semi-definiteness of $E$ and $R$, there exists a nonsingular matrix  $X$ such that
\begin{eqnarray*}
    XEX^T&=&\bmat{ & n_1      & n_2    & n_3 & n_4 \cr
                        n_1     & E_{11} & 0       & 0     & 0     \cr
                        n_2     & 0         & E_{22} & E_{23} & 0 \cr
                        n_3     & 0         & E_{23}^T & E_{33} & 0 \cr
                        n_4     & 0         & 0              & 0         & 0 \cr}, \\
XRX^T&=&\bmat{& n_1   & n_2    &  n_3 &    n_4    \cr
                        n_1 & 0 & 0 & 0 & 0 \cr
                        n_2 & 0 & R_{22} & R_{23} & 0 \cr
                        n_3 & 0 & R_{23}^T & R_{33} & 0 \cr
                        n_4 & 0 & 0 & 0 & 0 \cr}, \
       XB=\bmat{ &        \cr
                 n_1  & 0     \cr
                 n_2 & B_2  \cr
                 n_3 & 0      \cr
                 n_4 & 0     \cr},
\end{eqnarray*}
 where $E_{11}>0$, $\mat{cc} E_{22} & E_{23} \\ E_{23}^T & E_{33} \rix\geq 0$,
  $R_{33}>0$, $\mat{cc} R_{22} & R_{23} \\ R_{23}^T &R_{33} \rix\geq 0$, and $\rank(B_2)=n_2$.

  Set
\[ FX=\bmat{ & n_1 & n_2 & n_3 & n_4 \cr
                     & F_1 & F_2 & F_3 & F_4 \cr}.
                     \]
Then
\[
X\tilde R X^T=\mat{cccc} 0 & 0 & 0 & 0 \\
                       -{1\over 2} B_2F_1 & R_{22}-{1\over 2}(B_2F_2+F_2^TB_2^T) & R_{23}-{1\over 2}B_2F_3 & {1\over 2} B_2F_4 \\
                      0                            & R_{23}^T-{1\over 2}F_3^TB_2^T & R_{33} & 0 \\
                      0                            & -{1\over 2} F_4^TB_2^T & 0 & 0 \rix \geq 0
                    \]
gives that
\[ B_2F_1=0, \quad B_2F_4=0,
\]
which together with the condition (\ref{F2}) imply
\[ \rank \mat{cc} R_{22}-{1\over 2}(B_2F_2+F_2^TB_2^T) & R_{23}-{1\over 2}B_2F_3 \\
                          R_{23}^T-{1\over 2}F_3^TB_2^T & R_{33} \rix=n_2+n_3=\rank \mat{cc} R & B \rix, \]
and thus,
\[     \mat{cc} R_{22}-{1\over 2}(B_2F_2+F_2^TB_2^T) & R_{23}-{1\over 2}B_2F_3 \\
                          R_{23}^T-{1\over 2}F_3^TB_2^T & R_{33} \rix >0. \]
Consequently, there exists a nonsingular matrix $\cX$ such that
\[ \mat{ccc} I_{n_1} &      &     \\
                               & \cX &     \\
                               &       & I_{n_4} \rix X
  E X^T \mat{ccc} I_{n_1} &      &     \\
                               & \cX &     \\
                               &       & I_{n_4} \rix^T
                     =\bmat{ & n_1 & \tau_2 & n_2+n_3-\tau_2 & n_4 \cr
                         n_1    & E_{11} & 0 & 0 & 0 \cr
                        \tau_2 & 0  & \cE_{22} & 0 & 0 \cr
          n_2+n_3-\tau_2 & 0   & 0 & 0 & 0 \cr
                        n_4     & 0 & 0 & 0 & 0 \cr},  \]
\[   \mat{ccc} I_{n_1} &      &     \\
                               & \cX &     \\
                               &       & I_{n_4} \rix X
  \tilde R X^T \mat{ccc} I_{n_1} &      &     \\
                               & \cX &     \\
                               &       & I_{n_4} \rix^T
          =\bmat{ & n_1 & \tau_2 & n_2+n_3-\tau_2 & n_4 \cr
                         n_1    & 0 & 0 & 0 & 0 \cr
                        \tau_2 & 0  & \tilde R_{22} & \tilde R_{23} & 0 \cr
          n_2+n_3-\tau_2 & 0   & \tilde R_{23}^T & \tilde R_{33} & 0 \cr
                        n_4     & 0 & 0 & 0 & 0 \cr}, \]
\[   \mat{ccc} I_{n_1} &      &     \\
                               & \cX &     \\
                               &       & I_{n_4} \rix X B=\bmat{ &       \cr
              n_1     & 0   \cr
              \tau_2 & \cB_2 \cr
       n_2+n_3-\tau_2 & \cB_3 \cr
                   n_4 & 0 \cr}, \]
where
\[ E_{11}>0, \ \cE_{22}>0, \ \mat{cc} \tilde R_{22} & \tilde R_{23} \\ \tilde R_{23}^T & \tilde R_{33} \rix>0, \ \tilde R_{33}>0. \]
Set, with analogous partitioning,
\[    \mat{ccc} I_{n_1} &      &     \\
                               & \cX &     \\
                               &       & I_{n_4} \rix X J X^T \mat{ccc} I_{n_1} &      &     \\
                               & \cX &     \\
                               &       & I_{n_4} \rix^T  \]
\[                      =\bmat{ & n_1           & \tau_2       & n_2+n_3-\tau_2 & n_4 \cr
                               n_1 & J_{11}      & J_{12}      & J_{13}               & J_{14} \cr
                           \tau_2 & -J_{12}^T & J_{22}      & J_{23}              & J_{24} \cr
             n_2+n_3-\tau_2 & -J_{13}^T & -J_{23}^T & J_{33}              & J_{34} \cr
                               n_4 & -J_{14}^T & -J_{24}^T & -J_{34}^T        & J_{44} \cr}.
\]
Transforming and partitioning $\tilde J=J+{1\over 2}(BF-F^TB^T)$ as
\[
\mat{ccc} I_{n_1} &      &     \\
                               & \cX &     \\
                               &       & I_{n_4} \rix X\tilde J X^T \mat{ccc} I_{n_1} &      &     \\
                               & \cX &     \\
                               &       & I_{n_4} \rix^T   \]
\[      =\bmat{ & n_1           & \tau_2       & n_2+n_3-\tau_2 & n_4 \cr
                               n_1 & J_{11}      & J_{12}      & J_{13}               & J_{14} \cr
                           \tau_2 & -J_{12}^T & \tilde J_{22}      & \tilde J_{23}              & J_{24} \cr
             n_2+n_3-\tau_2 & -J_{13}^T & -\tilde J_{23}^T & \tilde J_{33}              & J_{34} \cr
                               n_4 & -J_{14}^T & -J_{24}^T & -J_{34}^T        & J_{44} \cr}.
\]
Then, using the condition (\ref{index-con}), it follows that
\[ \rank \mat{c} -J_{34}  \\  J_{44} \rix=n_4. \]
Hence, using Lemma \ref{lemma-1} (ii), we have  that
\[
\mat{cc} \tilde J_{33}-\tilde R_{33} & J_{34} \\ -J_{34}^T & J_{44} \rix
\]
is nonsingular and thus the pencil
\[  \mat{ccc} I_{n_1} &      &     \\
                               & \cX &     \\
                               &       & I_{n_4} \rix X (sE-(J-R+BF))
    X^T  \mat{ccc} I_{n_1} &      &     \\
                               & \cX &     \\
                               &       & I_{n_4} \rix^T  \]
is regular and of index at most one, or equivalently, the pencil $sE -(J-R+BF)$ is regular and of index at most one.
}
\eproof

Using these preliminary lemmas, in the next section we present our main results.
\section{Main Results}\label{sec:main}
In this section we present the answers to the two questions when a pH descriptor system can be made regular, index at most one and asymptotically stable as well as strictly passive. We present { constructive proofs  that can be implemented to give an explicit construction} of the desired feedback and we begin with Question 1.
%
\begin{theorem}\label{Th1}  Consider a pH descriptor system of the form \eqref{genphdae}. There exists  a matrix $F\in \Rmn$ {(and we present a construction of such an $F$ in the proof)} such that closed loop system \eqref{clo-sys1} is port-Hamiltonian,and the pencil $sE-( \tilde J-\tilde R)$ is regular, of index at most one,  and has all its finite eigenvalues in $\C_-$ if and only if
the following two conditions hold:
\begin{eqnarray}
  & &  \rank \mat{ccc} sE-(J-R) & (G-P)(S+N)^{(+)}\cR(S)  & (G-P)\cS_\infty(S+N) \rix=n, \nonumber \\
  & &   \mbox{\textrm for all}\ s\in \C_0,  \label{con1}
\end{eqnarray}
and
\begin{equation}\label{con1-2} \rank \mat{cccc} E & (J-R)\cS_\infty(E) & (G-P)(S+N)^{(+)}\cR(S) & (G-P)\cS_\infty(S+N)   \rix=n.
\end{equation}
\end{theorem}

\proof

System (\ref{genphdae}) is port-Hamiltonian, so, (\ref{lmi-1}) holds.
It follows from $S\geq 0$, see e.~g. \cite{AchAC21,BeaMXZ18}, that there exists an orthogonal matrix $U$ such that
\begin{equation}\label{D1}
    U^T(S+N)U=\bmat{ & m_1          & m_2     & m_3 \cr
    m_1   & D_{11}      & D_{12} & 0     \cr
    m_2   & -D_{12}^T & D_{22} & 0   \cr
    m_3   & 0               & 0          & 0   \cr}, \quad
U^TSU=\bmat{ & m_1  & m_2 & m_3 \cr
m_1    & S_{11}  & 0 & 0     \cr
m_2   &  0     & 0 & 0     \cr
m_3   &  0     & 0  & 0   \cr},
\end{equation}
where
\begin{equation}\label{D2}  D_{22}^T=-D_{22},  \quad
\rank \mat{cc} D_{11} & D_{12} \\
-D_{12}^T & D_{22} \rix=m_1+m_2, \quad \rank(S_{11})=m_1, \quad S_{11}>0.
\end{equation}
Partitioning the resulting matrices as
\begin{eqnarray*}
&& U=\bmat{ & m_1 & m_2 & m_3 \cr
& U_1 & U_2 & U_3 \cr},\
 PU=\bmat{ & m_1 & m_2 & m_3 \cr
        & P_1 & P_2 & P_3 \cr}, \\
&& (G-P)U \mat{ccc} D_{11} & D_{12} & 0 \\
        -D_{12}^T & D_{22} & 0 \\
    0           & 0   & I \rix^{-1}=\bmat{ & m_1 & m_2 & m_3 \cr              & B_1 & B_2 & B_3 \cr},
\end{eqnarray*}
we have
\begin{eqnarray*}
    U_3&=&\cS_\infty(S+N), \quad U_1=\cR(S), \\
    B_1&=&(G-P)(S+N)^{(+)}U_1=(G-P)(S+N)^{(+)}\cR(S),\\  B_3&=&(G-P)U_3=(G-P)\cS_\infty(S+N).
\end{eqnarray*}
Moreover,  by (\ref{lmi-1}) we have $P_2=0$ and $P_3=0$.

Partition also a transformed feedback matrix $F\in \Rmn$ as
\[
\mat{ccc} D_{11} & D_{12} & 0 \\ -D_{12}^T & D_{22} & 0 \\ 0 & 0 & I \rix U^T F=\bmat{ &          \cr
m_1   & F_1  \cr
m_2   & F_2  \cr
m_3   & F_3  \cr}.
\]

We first show the necessity of conditions~\eqref{con1} and~\eqref{con1-2}.
If the closed-loop  system (\ref{clo-sys1})  is a pH descriptor system, then
\[
\mat{cc} \tilde R & \tilde P \\ \tilde P^T & S \rix
    =\mat{cc} R-{1\over 2}[(G-P)F+F^T(G-P)^T] & P+{1\over 2} F^T(S+N)^T \\ P^T+{1\over 2}(S+N)F & S \rix \geq 0
\]
which in the partitioned form implies that
\[
\mat{cccc} R-{1\over 2}[B_1F_1+B_2F_2+B_3F_3+(B_1F_1+B_2F_2+B_3F_3)^T] & P_1+{1\over 2}F_1^T & {1\over 2}F_2^T & 0 \\
P_1^T+{1\over 2}F_1 & S_{11} & 0 & 0 \cr
{1\over 2}F_2 & 0 & 0 & 0 \cr                                    0 & 0 & 0 & 0 \rix\geq 0 .
\]
Thus  $F_2=0$ and
\[
\mat{cc} R-{1\over 2}[B_1F_1+B_3F_3+(B_1F_1+B_3F_3)^T] & P_1+{1\over 2}F_1^T  \\                               P_1^T+{1\over 2}F_1  & S_{11}    \rix\geq 0.
\]
It follows that
\[
(G-P)F=\mat{cc} B_1 & B_3 \rix \mat{cc} F_1 \\ F_3 \rix.
\]
Since the pencil  $sE-(J-R+(G-P)F)$ is regular and of index at most one, and also has all finite eigenvalues in $\C_-$, it follows that
{
\begin{eqnarray}
& & \rank \mat{ccc} sE-(J-R) & B_1  & B_3 \rix \nonumber \\ &=& \rank \mat{ccc} sE-(J-R+\mat{cc} B_1 & B_3 \rix \mat{c} F_1 \\ F_3 \rix) & B_1 & B_3 \rix \nonumber \\
       &=& \rank \mat{ccc} sE-(J-R+(G-P)F) & B_1 & B_3 \rix  \nonumber \\
       &=& n, \ \mbox{\textrm for all}\ s\in \C_0,  \label{con1-N}
\end{eqnarray}
}
and
\begin{equation}\label{con1-2N}
\rank \mat{cccc} E & (J-R)\cS_\infty(E) & B_1 & B_3   \rix=n,
\end{equation}
i.~e.,  the conditions (\ref{con1}) and  (\ref{con1-2}) are satisfied.

To show the sufficiency, assume that conditions (\ref{con1}) and  (\ref{con1-2}) hold, or equivalently that the conditions (\ref{con1-N}) and (\ref{con1-2N}) hold.
Since
\[
\mat{cc} R & P \\ P^T & S \rix\geq 0,
\]
also
\[
\mat{cc} I & 0 \\ 0 & U \rix^T \mat{cc} R & P \\ P^T & S \rix \mat{cc} I & 0 \\ 0 & U \rix \geq 0,
\]
which gives
\[
\mat{cc} R & P_1 \\ P_1^T & S_{11} \rix \geq 0,
\]
and thus
\begin{eqnarray*}
& & \mat{cc} R+B_1P_1^T+P_1B_1^T+B_1S_{11}B_1^T & P_1+B_1S_{11} \\
                    P_1^T+S_{11} B_1^T & S_{11} \rix      \\
&=&   \mat{cc} I & B_1 \\ 0 & I \rix \mat{cc} R & P_1 \\ P_1^T & S_{11} \rix  \mat{cc} I & B_1 \\ 0 & I \rix^T \geq 0.
\end{eqnarray*}
We obtain
\begin{equation}\label{con1-1}  R+B_1P_1^T+P_1B_1^T+B_1S_{11}B_1^T \geq 0,
\end{equation}
and can construct a feedback matrix via
\[
F_2=0, \quad F_1=-2(B_1S_{11}+P_1)^T.
\]
Then
\begin{eqnarray*}
\tilde R&=&R-{1\over 2}(G-P)F-{1\over 2}F^T(G-P)^T\\
&=&R-{1\over 2}B_3F_3-{1\over 2} F_3^TB_3^T+(B_1P_1^T+P_1B_1^T+2B_1S_{11}B_1^T).
\end{eqnarray*}
Since $R\geq 0$, we can determine a nonsingular matrix $Z$ and orthogonal matrices $V_1$ and $V_3$ such that
\begin{eqnarray*}
ZB_3V_3&=&\bmat{ & \mu_1 & m_3-\mu_1 \cr
                     \mu_1 & I  & 0      \cr
                    \mu_2  & 0 & 0 \cr
                    \mu_3 & 0 & 0 \cr
                    \mu_4 & 0 & 0 \cr}, \quad
  ZB_1S_{11}^{1/2}V_1=\bmat{ & \mu_2 & m_1-\mu_2 \cr
                \mu_1  & 0 & B_{12} \cr
                \mu_2  & I  & 0   \cr
                \mu_3  & 0 & 0 \cr
                \mu_4 & 0 & 0 \cr}, \\
                ZRZ^T&=&\bmat{ & \mu_1 & \mu_2 & \mu_3 & \mu_4 \cr
           \mu_1 & R_{11} & R_{12} & 0 & 0 \cr
           \mu_2 & R_{12}^T & R_{22} & 0 & 0 \cr
           \mu_3 & 0 & 0 & I & 0 \cr
           \mu_4 & 0 & 0 & 0 & 0 \cr},
\end{eqnarray*}
where due to the congruence transformation
\[
\mat{cc} R_{11} & R_{12} \\ R_{12}^T & R_{22} \rix\geq 0.
\]
Partition the transformed part $F_3$ of the feedback matrix as
\[
V_3^TF_3Z^{T}=\bmat{ & \mu_1 & \mu_2 & \mu_3 & \mu_4 \cr
\mu_1 & F_{31} & F_{32} & F_{33} & F_{34} \cr
    m_3-\mu_1  & 0 & 0 & 0 & 0 \cr},
\]
with $F_{31}$, $F_{32}$, $F_{33}$  and $F_{34}$ to be determined and set
\[
V_1^TS_{11}^{-1/2} P_1^TZ^T=\bmat{ & \mu_1 & \mu_2 & \mu_3 & \mu_4 \cr
                        \mu_2    & P_{11} & P_{12} & P_{13} & P_{14} \cr
                     m_1-\mu_2& \hat P_{11} & \hat P_{12} & \hat P_{13} & \hat P_{14} \cr}.
\]
Then it follows from (\ref{con1-1}) that
\begin{eqnarray*}
  &&  Z(R+B_1P_1^T+P_1B_1^T+B_1S_{11}B_1^T)Z^T \\
& = & \mat{cccc} R_{11}+B_{12}\hat P_{11}+\hat P_{11}^TB_{12}^T+B_{12}B_{12}^T & R_{12}+B_{12}\hat P_{12}+P_{11}^T & B_{12}\hat P_{13} & B_{12}\hat P_{14} \\
                      R_{12}^T+\hat P_{12}^TB_{12}^T+P_{11} & R_{22}+P_{12}+P_{12}^T+I  & P_{13} & P_{14} \\
                      \hat P_{13}^TB_{12}^T & P_{13}^T & I & 0 \\
                      \hat P_{14}^TB_{12}^T & P_{14}^T & 0 & 0 \rix  \\
  &\geq&  0,
\end{eqnarray*}
which implies that
\[
B_{12}\hat P_{14}=0, \quad P_{14}=0, \]
and
\[ \mat{cc} R_{22}+P_{12}+P_{12}^T +I & P_{13} \\ P_{13}^T & I \rix \geq 0,
\]
consequently,
\[ R_{22}+P_{12}+P_{12}^T +I-P_{13}P_{13}^T \geq 0, \]
\begin{equation} \label{add1}
    R_{22}+P_{12}+P_{12}^T +2I-P_{13}P_{13}^T=  (R_{22}+P_{12}+P_{12}^T +I-P_{13}P_{13}^T) +I>0,
\end{equation}
and furthermore,
\begin{equation}\label{A1}  \mat{cc} R_{22}+P_{12}+P_{12}^T +2I  & P_{13} \\ P_{13}^T & I \rix>0, \quad
     \rank  \mat{cc} R_{22}+P_{12}+P_{12}^T +2I  & P_{13} \\ P_{13}^T & I \rix=\mu_2+\mu_3.
\end{equation}
A simple calculation yields that
\begin{eqnarray*}
Z \tilde R  Z^T
&=& \mat{cccc}
  R_{11}+B_{12}\hat P_{11}+\hat P_{11}^TB_{12}^T+2B_{12}B_{12}^T & R_{12}+B_{12}\hat P_{12}+P_{11}^T & B_{12}\hat P_{13} & 0 \\
  R_{12}^T+\hat P_{12}^TB_{12}^T+P_{11} & R_{22}+P_{12}+P_{12}^T+2I & P_{13} & 0 \\
  \hat P_{13}^TB_{12}^T & P_{13}^T & I & 0 \\
  0 & 0 & 0 & 0 \rix \\
& & -{1\over 2} \mat{cccc}
     F_{31}+F_{31}^T  & F_{32} & F_{33} & F_{34} \\
     F_{32}^T & 0  & 0 & 0 \\
     F_{33}^T & 0  & 0 & 0 \\
     F_{34}^T & 0 & 0 & 0 \rix,
\end{eqnarray*}
and we can choose
\[   F_{32}=2(R_{12}+B_{12}\hat P_{12}+P_{11}^T), \quad  F_{33}=2B_{12}\hat P_{13}, \quad  F_{34}=0, \]
and any $F_{31}$  such that
\[
R_{11}-{1\over 2}(F_{31}+F_{31}^T)+B_{12}\hat P_{11}+\hat P_{11}^TB_{12}^T+B_{12}B_{12}^T>0.
\]
We then have
\[ \mat{cc} Z & 0 \\ 0 & I \rix
\mat{cccc} I & B_1 & 0 & 0 \\                          0 & I & 0 & 0 \\                                        0 & 0 & I & 0 \\
0 & 0 & 0 & I \rix \mat{cc} I & 0 \\ 0 & U^T \rix
  \mat{cc} \tilde R & \tilde P \\ \tilde P^T & S \rix (\mat{cc} Z & 0 \\ 0 & I \rix \mat{cccc} I & B_1 & 0 & 0 \\                                              0 & I & 0 & 0 \\                                     0 & 0 & I & 0 \\                                     0 & 0 & 0 & I \rix \mat{cc} I & 0 \\ 0 & U^T \rix)^T \]
%
%
\[  =\mat{ccccccc}  R_{11}-{1\over 2}(F_{31}+F_{31}^T)+B_{12}\hat P_{11}+\hat P_{11}^TB_{12}^T+B_{12}B_{12}^T & 0 & 0 & 0 & 0 & 0 & 0\\
                         0 & R_{22}+P_{12}+P_{12}^T+I  & P_{13} & 0 & 0 & 0 & 0 \\
                         0 & P_{13}^T & I & 0 & 0 & 0 & 0 \\
                         0 & 0 & 0 & 0 & 0  & 0 & 0  \\
                         0 & 0 & 0 & 0 & S_{11} & 0 & 0 \\
                         0 & 0 & 0 & 0 & 0                  & 0 & 0 \\
                         0 & 0 & 0 & 0 & 0                  & 0 & 0 \rix   \]
\[ \geq 0,  \]
and thus
\[
\mat{cc} \tilde R & \tilde P \\ \tilde P^T & S \rix \geq 0,
\]
i.~e. the  closed-loop system (\ref{clo-sys1})  is a  pH descriptor system.  Moreover,
\begin{eqnarray}
& & Z\tilde R Z^T \nonumber \\
&=& Z (R-{1\over 2} \mat{cc} B_1 & B_3 \rix \mat{c} F_1 \\ F_3 \rix-{1\over 2} (\mat{cc} B_1 & B_3 \rix \mat{c} F_1 \\ F_3 \rix)^T ) Z^T    \nonumber \\
&=& \mat{cccc}  R_{11}-{1\over 2}(F_{31}+F_{31}^T)+B_{12}\hat P_{11}+\hat P_{11}^TB_{12}^T+2B_{12}B_{12}^T & 0 & 0 & 0  \\
                         0 & R_{22}+P_{12}+P_{12}^T+2I  & P_{13} & 0  \\
                         0 & P_{13}^T & I & 0 \\
                         0 & 0 & 0 & 0 \rix     \nonumber \\
&\geq& 0, \label{R}
\end{eqnarray}
and
\begin{eqnarray*}
\rank(\tilde R) &=&  \rank(R_{11}-{1\over 2}(F_{31}+F_{31}^T)+B_{12}\hat P_{11}+\hat P_{11}^TB_{12}^T+2B_{12}B_{12}^T) \\
& &  +
\rank \mat{cc}  R_{22}+P_{12}+P_{12}^T+2I  & P_{13} \\   P_{13}^T & I  \rix    \\
&=& \mu_1+\mu_2+\mu_3                                                      \\
&=&   \rank \mat{ccc} R & B_1 & B_3 \rix.
\end{eqnarray*}
Hence, it follows from Lemmas \ref{lemma-2} and  \ref{lemma-3} that the pencil $sE-(J-R+(G-P)F)$ is regular, of index at most one, and has all eigenvalues in the open left half plane.  \eproof

Question 2. is answered in the following theorem.
\begin{theorem}\label{Th2} There exists a matrix $F$ { (and we present a construction of such $F$ in the proof)} such that the  closed-loop system (\ref{clo-sys1})  is strictly passive if and only if $S>0$ and
\begin{equation}\label{con2}
R+{1\over 2} (G-P)(S+N)^{-1}(G+P)^T+{1\over 2} (G+P) (S+N)^{-T}(G-P)^T>0.
\end{equation}
\end{theorem}
\proof  By setting $B=G-P$, $D=S+N$, {and using the notation as in Theorem~\eqref{Th1}}, we have the following list of equivalences.
\begin{eqnarray*}
& &  \mat{cc} \tilde R & \tilde P \\ \tilde P^T & S \rix>0 \\
&\Leftrightarrow& \mat{cc} R-{1\over 2}BF-{1\over 2}F^TB^T & P+{1\over 2}F^TD^T \\ P^T+{1\over 2}DF & S \rix>0 \\
&\Leftrightarrow& \mat{cc} I & BD^{-1} \\ 0 & I \rix  \mat{cc} R-{1\over 2}BF-{1\over 2}F^TB^T & P+{1\over 2}F^TD^T \\ P^T+{1\over 2}DF & S \rix \mat{cc} I & BD^{-1} \\ 0 & I \rix^T>0 \\
&\Leftrightarrow& \mat{cc} R+BD^{-1}P^T+PD^{-T}B^T+BD^{-1}SD^{-T}B^T & P+BD^{-1}S+{1\over 2}F^TD^T \\
P^T+SD^{-T}B^T +{1\over 2}DF & S \rix  >0 \\
&\Leftrightarrow& S>0, \mbox{\textrm and }\ R+BD^{-1}P^T+PD^{-T}B^T+BD^{-1}SD^{-T}B^T >0,
\end{eqnarray*}
i.~e., the condition (\ref{con2}) holds (since $S={1\over 2} (D+D^T)$). Moreover, under the condition (\ref{con2}),
\[ F = -2(S+N)^{-1}[P^T+S(S+N)^{-T}B^T]=-(S+N)^{-1} (G+P)^T-(S+N)^{-T}(G-P)^T \]
is a feedback matrix that achieves strict passivity of  the closed loop system. \eproof

\section{Conclusions and Future Work}
\label{sec:conclusions}
{
We have presented a complete characterization when a pH descriptor system (\ref{genphdae}) can be made asymptotically stable and strictly passive, as well as regular and of index at most one via state feedback while  preserving the pH structure.

A complete characterization  when a pH descriptor system (\ref{genphdae}) with $D=0$  can be made regular, of index at most one and asymptotically stable by output feedback while preserving the pH structure has been presented in \cite{ChuM24a},  while the question how to achieve  a regular, index at most one and  asymptotically stable or strictly passive system via output feedback for a pH descriptor system (\ref{genphdae}) with  feedthrough term $D\not =0$ is still an open problem.

An interesting question that also remains to be answered is to use the discussed state feedback to obtain a representation that maximizes the distance to instability, respectively non-passivity. For both questions only partial solutions are available, see \cite{MehU23} and the references therein.}


\end{document}